\documentclass[12pt]{amsart}

\usepackage{amsmath, amssymb, amsthm}
\usepackage{graphicx}
\usepackage{amsfonts, stmaryrd, mathtools, relsize, tikz-cd, mathrsfs, extarrows, mathrsfs}
\usepackage{cleveref}
\usepackage{comment}
\usepackage{ytableau}
\usepackage{float}
\usepackage[margin=1in]{geometry}
\usepackage{adjustbox}

\usepackage[numbers]{natbib}

\usepackage{wrapfig}

\newcommand{\Pp}{\mathbb{P}}

\newcommand{\id}{\mbox{id}}
\newcommand{\coker}{\mbox{coker}}

\newcommand{\Sym}{\mbox{Sym}}
\newcommand{\Oo}{\mathbb{O}}
\newcommand{\Ss}{\mathbb{S}}

\newtheorem{theorem}{Theorem}[section]
\newtheorem{lemma}[theorem]{Lemma}
\newtheorem{proposition}[theorem]{Proposition}
\newtheorem{corollary}[theorem]{Corollary}

\theoremstyle{definition}
\newtheorem{example}[theorem]{Example}

\theoremstyle{remark}
\newtheorem{remark}[theorem]{Remark}
\newtheorem*{remark*}{Remark}

\numberwithin{equation}{section}
\begin{document}
\bibliographystyle{alpha}
\title[Self-Duality of Generalized Eagon-Northcott Complexes]{Hermite Reciprocity and Self-Duality of Generalized Eagon-Northcott Complexes}

\author{Ethan Reed}
\address{Department of Mathematics, University of Notre Dame, 255 Hurley, Notre Dame, IN 46556}
\email{ereed4@nd.edu}

\date{\today}

\subjclass[2010]{Primary 13D02, 13D07, 14F06, 15A69, 18G40, 20G05}

\keywords{}

\begin{abstract}
    Previous examples of self-duality for generalized Eagon-Northcott complexes were given by computing the divisor class group for Hankel determinantal rings. We prove a new case of self-duality of generalized Eagon-Northcott complexes with input being a map defining a Koszul module with nice properties. This choice of Koszul module can be specialized to the Weyman module, which was used in a proof of the generic version of Green's conjecture. In this case, the proof uses a version of Hermite Reciprocity not previously defined in the literature. 
\end{abstract}

\maketitle

\section{Introduction}
We begin with a discussion of Hermite Reciprocity, a classical result in the invariant theory of binary forms, i.e. the study of $SL_2(\mathbb{C})$-representations. If $U = k^2$ is the standard representation of $SL_2(k)$, then Hermite Reciprocity is the statement that the following representations are isomorphic:
\[\Sym^{m}(\Sym^n U) \cong \bigwedge^m(\Sym^{m+n-1}U) \cong \Sym^n(\Sym^m U)\]
Often given as an exercise in an introductory course in representation theory, this statement can be proven by enumerating the indecompasable factors occurring in each term \cite{FultonHarris}[Exercise 11.35]. A more subtle question is to ask if instead one can write down natural maps, which give these isomorphisms. In \cite{RaicuSam}, Raicu and Sam discuss instances in which Hermite Reciprocity isomorphisms arose in the study of Green's Conjecture in \cite{Voisin}\cite{RaicuSam0}, and  \cite{Greens}, and new instances in connection with the Hankel determinantal ring divisor class group, Schwarzenberger bundles, and the secant varieties to rational normal curves. In each case, the abstract identification of Hermite Reciprocity was not sufficient, but rather a particular identification via natural maps was needed. The authors further show that in each case the natural maps that appeared were all equivalent. In the proof of the main result of this paper, Hermite Reciprocity again arises, but this time it is not equivalent to prior formulations. As the definitions of the maps giving the Hermite isomorphisms is rather technical, we do not give their exact formulation in the introduction. A comparison of this new special Hermite isomorphism with the formulations discussed in \cite{RaicuSam} is done in \Cref{specHerm}.

The main result of this paper concerns the complex of free modules $\Sym^{f-g+1}(\varphi)$ constructed from a map of free modules 
\[\varphi: S^f = F \to G = S^g,\]
and is given by
\[\wedge^{f-g+1}F \to \cdots \to F \otimes \mbox{Sym}^{f-g}(G) \to \mbox{Sym}^{f-g+1}(G) \to 0\]
These complexes are examples of generalized Eagon-Northcott complexes, which are a family of complexes first introduced in \cite{EisenbudBuchsbaum} (The case above being $\mathcal{C}_{f-g+1}$ in \cite{Eisenbud}[Appendix 2.6]). Generalized Eagon-Northcott complexes show up in a number of contexts in commutative algebra and algebraic geometry with additional generalizations (for example: \cite{EN1}, \cite{EN2}, \cite{EN3}, \cite{EN4}, \cite{EN5}). Even though the ranks of the free modules in the complex $\mbox{Sym}^{f-g+1}(\varphi)$ are symmetric, this complex is not necessarily self-dual. In particular, in the universal setting these complexes are not self-dual for $g > 1$ (see \Cref{notdual}).

Despite these complexes not being self-dual in the universal setting, these complexes can be self-dual for certain specializations of $\varphi$. For example consider the case coming from the representation theory of binary forms when $U$ is a two dimensional $\mathbb{C}$-vector space, $S = \mbox{Sym}(\mbox{Sym}^{d+b}(U)) = \mathbb{C}[x_0 \cdots x_{d+b}]$, and $\varphi$ is the linear map,
\[\varphi: \mbox{Sym}^b(U) \otimes S \to \mbox{Sym}^d(U) \otimes S\]
defined by the map of $SL_2(\mathbb{C}) = SL(U)$-representations
\[\mbox{Sym}^b(U) \to \mbox{Sym}^d(U)\otimes \mbox{Sym}^{d+b}U\]
arising from an inclusion of a factor of highest weight.
After a choice of basis, $\varphi$ can be written as the Hankel matrix
\[\varphi = \begin{bmatrix} 
x_0 & \cdots & x_d\\
x_1 & \cdots & x_{d+1}\\
\vdots & & \vdots\\
x_b & \cdots & x_{d+b}
\end{bmatrix}\]

In this case, it is known that the complex $\mbox{Sym}^{b-1}(\varphi)$ and its dual are minimal free resolutions for modules that are isomorphic as noted in Remark 3.7 of \cite{Hankel}. More is true in this case, as the isomorphism of complexes can be made to be $SL_2(\mathbb{C})$-equivariant as shown in \cite{RaicuSam}. This is again done by showing that the modules that are minimally resolved by these complexes are isomorphic, but now respecting the $SL_2(\mathbb{C})$-action. Critical to this proof is the Hermite Reciprocity isomorphism that has shown up in a number of contexts as mentioned above.

The following is a special case of the main result of this paper.
\begin{theorem}\label{special1}
Let $U$ be the standard representation of $\mbox{SL}_{2}(\mathbb{C})$, $b \geq 2$, and $S = \mbox{Sym}^{\bullet}(\mbox{Sym}^bU)$. Let $\varphi$ be the linear map
$$\varphi:\mbox{Sym}^{2b-2}U \otimes S \to \mbox{Sym}^bU \otimes S$$
defined by the inclusion of representations of second highest weight
\[\mbox{Sym}^{2b-2}U \subset \mbox{Sym}^b U \otimes \mbox{Sym}^bU.\]
Then the complex $\mbox{Sym}^{b-1}(\varphi)$ is self-dual.
\end{theorem}
 If we let $V_0 = \Sym^{b} U$, then 
 \[V_1:= \Sym^{2b-2}U \subset \wedge^2(\Sym^bU) =  \wedge^2 V_0\]
 by the Clebsch-Gordan rule for wedge products \cite[Exercise 11.30]{FultonHarris}.
 Further, the map 
 \[\Sym^{2b-2} U \otimes S \to \Sym^{b} U \otimes S\]
 fits into a sequence 
 \[\Sym^{2b-2}\otimes S \to \Sym^{b} U \otimes S \to S,\]
 where the middle homology of this sequence is the Weyman module used in a proof of the generic version of Green's conjecture in \cite{Greens}. To prove that the Weyman module satisfied a certain vanishing statement, a vanishing statement was instead proven for a more general construction of Koszul modules, where $\Sym^{b-2}U$ was replaced with other subspaces $V_1 \subset \wedge^2 V_0$. This is similar to what we do here to obtain \Cref{special1} as a special case of \Cref{main1} 
\begin{theorem}\label{main1}
 Let $b \geq 1$, $V_0$ be a $b+1$-dimensional complex vector space, $S = Sym(V_0) = \mathbb{C}[x_0, ... , x_b]$, and $m = (x_0, ... , x_b)$. Recall that $\mathbb{C} = S/m$ is resolved by the Koszul complex, which ends in a sequence 
\[\wedge^2V_0 \otimes S \xrightarrow[]{\varphi} V_0 \otimes S \to S \to 0\]
Let $V_1 \subset \wedge^2V_0$ of dimension $2b-1$ such that in the restricted sequence
\[V_1 \otimes S \xrightarrow[]{\varphi|_{V_1 \otimes S}} V_0 \otimes S \to S \to 0\]
the middle homology is of finite length.

Then the complex $\mbox{Sym}^{b-1}(\varphi|_{V_1\otimes S})$ is self-dual.
\end{theorem}
Note that in this case the statement that the sequence has middle homology of finite length is equivalent to the corresponding sequence of sheaves being exact. This is the formulation that is used in most of the proofs. Unlike in the case of the Hankel matrix, the complexes in \Cref{main1} are no longer exact, so a different technique is needed to prove the isomorphism. As these complexes are linear, they correspond to graded modules over the exterior algebra $\wedge^{\bullet} (V_0^*)$ by the Bernstein-Gel'fand-Gel'fand correspondence \cite{syzygies}[Chapter 7].  It is then sufficient to instead show that thse modules are isomorphic. Critical to this proof is a naturally defined isomorphism
\[\wedge^{b-1}(V_1) \xrightarrow[]{\cong} \Sym^{b-1}(V_0),\]
which in the case $V_0 \Sym^b U$ and $V_1 = \Sym^{2b-2}U$ becomes a Hermite Reciprocity isomorphism
\[\wedge^{b-1}(\Sym^{2b-2} U) \otimes S \to \Sym^{b-1}(\Sym^b U)\]
not previously defined in the literature as mentioned above. 

The paper is outlined as follows. In \Cref{prelimns}, we construct the family of generalized Eagon-Northcott complexes and state some of the interesting properties of these complexes. In \Cref{embedding}, we embed the complex $\Sym^{b-1}(\varphi|_{V_1 \otimes S})$ into the minimal free resolution of a power of the maximal ideal of $S$. In \cref{BGG}, it is shown that the graded module over $\wedge^{\bullet}(V_0)$ that corresponds to $\Sym^{b-1}(\varphi|_{V_1 \otimes S})$ is generated in lowest degree. In \Cref{proofsec}, the proof of \Cref{main1} is completed by giving a map from the dual of $\Sym^{b-1}(\varphi|_{V_1\otimes S})$ to the resolution in \Cref{embedding} and containing the image of the map defined there. By \Cref{BGG}, this containment is reduced to only needing to be checked in a single degree, and can be verified using a Hermite Reciprocity isomorphism in the case arising from the invariant theory of binary forms.

\section{Preliminaries}\label{prelimns}
In this section, we will define the family of Generalized Eagon-Northcott complexes using Schur complexes of two term complexes. This is in part to motivate \Cref{main1}, but also because these particular Schur complexes will be used again. For an overview of Schur complexes see \cite{Weyman}[Section 2.4], the construction here will be much less general. Let \[\varphi: F \to G\]
be a map of free modules over a commutative ring $R$, which we interpret as a chain complex with $G$ in degree $0$ and $F$ in degree $1$. We will define the Schur complexes, $\Sym^i(\varphi)$ and $\wedge^i(\varphi)$, using symmetric, exterior, and divided powers (denoted $\mbox{Sym}^j$, $\wedge^j$, $D^j$) of $F$ and $G$. 

The Schur complex $\Sym^i(\varphi)$ is the complex
\[\wedge^i F \to \cdots \to \wedge^2 F \otimes \Sym^{i-2}G \to F \otimes \Sym^{i-1}G\to \Sym^i G \to 0\]
where the differentials are given by the compositions 
\[\wedge^jF \otimes \Sym^{i-j}G \xrightarrow{\Delta \otimes \mbox{id}} \wedge^{j-1}F \otimes F \otimes \Sym^{i-j} \xrightarrow{\id \otimes \varphi \otimes \mbox{id}} \wedge^{j-1} F \otimes G \otimes \Sym^{i-j} G\]
\[ \xrightarrow{\id \otimes m} \wedge^{j-1} F \otimes \Sym^{i-j+1}\]
where $\Delta$ is the comultiplication map on the exterior algebra and $m$ is the multiplication map on the symmetric algebra. We now define the complex $\wedge^i(\varphi)$, which is defined similarly and is given by
\[D^i F \to \cdots \to D^2 F \otimes \wedge^{i-2} G\to F \otimes \wedge^{i-1}G\to \wedge^i G\]
The differentials of the complex are given by
\[D^j(F) \otimes \wedge^{i-j}G \xrightarrow{\Delta \otimes \mbox{id}} D^{j-1}F \otimes F \otimes \wedge^{i-j} \xrightarrow{\id \otimes \varphi \otimes \mbox{id}} D^{j-1} F \otimes G \otimes \wedge^{i-j} G\]
\[ \xrightarrow{\id \otimes m} D^{j-1} F \otimes \wedge^{i-j+1}\]
where $\Delta$ is now comultiplication on the divided power algebra and $m$ is multiplication on the exterior algebra. Up to a shift in homological degree, we have an identification independent of characteristic,
\begin{equation}\label{duality1}
\mbox{Sym}^i(\varphi) \cong (\wedge^i(\varphi^*))^*
\end{equation}
With these two term Schur complexes defined, we are ready to construct the generalized Eagon-Northcott complexes, a family of complexes 
\[C_i(\varphi)\]
associated to the map $\varphi$. For an algebraic overview of these complexes see \cite{Eisenbud}[Appendix 2.6] or for a more geometric approach \cite{Positivity}[Appendix B]. Let $f$ be the rank of $F$ and $g$ be the rank of $G$. Then for each $i \in \mathbb{Z}$, we construct the complex $C_{i}(\varphi)$ from the complexes $\Sym^i(\varphi)$ and $\wedge^{f-g-i}(\varphi*)$. Note first that by identifying $\wedge^f F \xrightarrow[]{\cong} R$, we can make the identification $(\wedge^{j} F)^* \cong \wedge^{f-j} F$. Using this identification, we can reformulate the complex $\wedge^{f-g-i}(\varphi^*)$ as
\begin{equation}\label{wedgedual}
    D^{f-g-i}(G^*)\to \cdots \to \wedge^{g+i+2}F \otimes D^2(G^*) \to \wedge^{g+i+1}F\otimes G^* \to \wedge^{g+i}F\end{equation}
where the differentials can now be written as
\begin{equation}\label{wedgedualdiff}D^j(G^*) \otimes \wedge^{g+i+j}F \xrightarrow[]{\Delta \otimes \Delta} D^{j-1}(G^*) \otimes G^* \otimes F \otimes \wedge^{g+i+j-1}F \xrightarrow[]{\id \otimes \varphi' \otimes \id}D^{j-1}G \otimes \wedge^{g+i+j-1}F\end{equation}
Here $\varphi'$ is the composition 
\[G^* \otimes F \to \xrightarrow[]{\id \otimes \varphi} G^* \otimes G \to R\]
This reformulation can be seen from the fact that $\Sym^{i}(\varphi)$ is a linear strand of the Koszul complex, the identification \Cref{duality1}, and the self-duality of the Koszul complex.

The terms of the complex $C_i(\varphi)$ are given by the terms of $\Sym^i(\varphi)$ and $\wedge^{f-g-i}(\varphi^*)$. More specifically, let 
\[C_i(\varphi)_{j} = \begin{cases}
\Sym^{i}(\varphi)_{j} & j \leq i\\
\wedge^{f-g-i}(\varphi)_{j-i-1} & j \geq i +1
\end{cases}\]
The differentials for $C_i(\varphi)$ are also given by the differentials of $\Sym^i(\varphi)$ and $\wedge^{f-g-i}(\varphi^*)$ except in homological degree $i+1$, in which case the differential is given by the composition
\[\wedge^{g+i} F \xrightarrow[]{\Delta} \wedge^gF \otimes \wedge^iF \xrightarrow[]{\wedge^g \varphi\otimes \id} \wedge^g G \otimes \wedge^i F \to \wedge^i F\]
This differential is typically referred to as a "splice" map. 

For an example of a complex of this family, consider the case of $C_{1}(\varphi)$, which is the Buchsbaum-Rim complex,
\[D^{f-g-1} (G^*)\otimes \wedge^{f} F \to D^{f-g-2}(G^*) \otimes \wedge^{f-1} F \to \cdots\]
\[\to \wedge^{g+1} F \to F \to G \to 0\]

Using \Cref{duality1} and a compatibility of the splice maps, we obtain another duality statement (up to a shift in homological degrees)
\begin{equation}\label{duality2} C_i(\varphi) \cong C_{f-g-i}(\varphi)^* \end{equation}
A general question about these complexes is when are these complexes self-dual, i.e. when is 
\begin{equation}\label{duality3} C_i(\varphi) \cong C_{f-g-i}(\varphi) \end{equation}
In each homological degree of $C_i(\varphi)$, the ranks of the free modules are given by products of binomial coefficients in terms of $f$ and $g$. By comparing these ranks for $g \geq 2$, \Cref{duality3} is only possible in the case when $i = f-g+1$ (or by duality $i = -1$). 
\begin{remark}\label{notdual}
    In the generic case when 
\[R = ZZ[x_{i,j}]_{\substack{1 \leq i \leq g, \\ 1 \leq j\leq f}}\]
and $\varphi$ is the map
\[\varphi: F \xrightarrow[]{(x_{i,j})} G,\]
these complexes are not self-dual, that is
\[C_{f-g+1}(\varphi) \ncong C_{-1}(\varphi)\]
In this case, both of these complexes are minimal free resolutions of non-isomorphic rank 1 Cohen Macaulay modules over the ring $R$ quotiented by the ideal of maximal minors of $[x_{i,j}]$. One way to see that these modules are not isomorphic is that they represent different elements of the divisor class group over this determinantal ring (see \cite{BrunsVetter}[Corollary 8.4]). For this reason, the statement of \Cref{main1} is in some sense unexpected.\end{remark}

Note also that if $i \geq f-g+1$, then $C_i(\varphi)$ is just the complex $\Sym^i(\varphi)$ and if $i \leq -1$, $C_{i}(\varphi)$ is just the complex $\wedge^{f-g+1}(\varphi)$, so in these cases we can adopt this notation to refer to the complexes $C_i(\varphi)$.

\section{Embedding $\Sym^{b-1} \varphi$ into an exact complex}\label{embedding}
Throughout the rest of the paper fix the notation that $V_0 = \mathbb{C}^{b+1}, S = \Sym^{\bullet} (V_0) \cong \mathbb{C}[z_0...z_b]$, and  $\varphi$ is the Koszul differential $\wedge^2 V_0 \otimes S \to V_0 \otimes S$. In this section, we define an $S$-degree preserving chain map from the complex $\mbox{Sym}^{b-1}(\varphi |_{V_1 \otimes S})$ to the minimial free resolution of $m^{b-1}$ for any choice of $V_1 \subset \wedge^2 V_0$. Under additional hypotheses on $V_1 \subset \wedge^2 V_0$, this map will be an embedding. 

$m^{b-1}$ has $gl(V_0)$-equavariant resolution given by 
$$\Ss(b-1, 1^b)V_0\otimes S \to \Ss(b-1, 1^{b-1})V_0 \otimes S$$
$$ \to ... \to \Ss(b-1, 1^2)V_0 \otimes S \to \Ss(b-1, 1)V_0 \otimes S \to\Ss(b-1) V_0 \otimes S \to 0$$
Denote this resolution by $C_{\bullet}$. To define a chain map preserving the degrees of $S$, $f(V_1)_{\bullet}:\Sym^{b-1} (\varphi|_{V_1 \otimes S}) \to C_{\bullet}$, it is sufficient to define this map on generators in each homological degree, which is requires maps 
\[f_i: \wedge^i V_1 \otimes \Sym^{b-1-i} V_0 \to \Ss(b-1, 1^i) V_0 \mbox{ for } 0 \leq i \leq b-1.\]
Note that $C_{\bullet}$ has one more non-zero term than $\Sym^{b-1} \varphi$, so we also have $f_{b}: 0 \to \Ss(b-1, 1^b) V_0$. Recall that 
\begin{equation}\label{Schurdef}\Ss(b-1, 1^i) V_0 = \ker(\wedge^i V_0 \otimes \Sym^{b-1} V_0 \xlongrightarrow{\delta} \wedge^{i-1} V_0 \otimes \Sym^{b} V_0),\end{equation}
where the map $\delta$ is the differential in the Koszul complex, i.e. is given by comultiplication on the exterior algebra followed by multiplication on the symmetric algebra. To define these maps, we first define a more general collection of maps $\psi_{i,j}(V_1)$ as follows. In the case of $V_1 = \wedge^2 V_0$ define
\[\psi_{i,j}(\wedge^2V_0): \wedge^{i}(\wedge^2 V_0) \otimes \Sym^{j}V_0 \to \mathbb{S}_{i+j, 1^j}V_0\]
as an $SL(V_0)$-equavariant projection onto an irreducible factor. Further define $\psi_{i,j}(V_1)$ to be the restriction of $\psi_{i,j}(\wedge^2 V_0)$ to 
\[\wedge^i(V_1) \otimes \Sym^j(V_0).\]
Using the following lemma, we obtain a more convenient for our purposes description of the maps $\psi_{i,j}(V_1)$.
\begin{lemma}\label{psidef}
 Using the inlcusion $V_1 \subset \wedge^2 V_0 \subset V_0 \otimes V_0$, the map $\psi_{i,j}(V_1)$ factors as in the diagram
 \[\begin{tikzcd}
 \wedge^iV_1 \otimes \Sym^jV_0 \arrow[d] \arrow[dddr, bend left = 30, "\psi_{i,j}(V_1)"]\\
 \wedge^i(V_0 \otimes V_0) \otimes \Sym^jV_0 \arrow[d]\\
 \wedge^i V_0 \otimes \Sym^iV_0 \otimes \Sym^j V_0 \arrow[d]\\
 \wedge^i V_0 \otimes \Sym^{i+j}V_0 & \arrow[l] \mathbb{S}_{i + j, 1^j} V_0
 \end{tikzcd}\]
\end{lemma}
\begin{proof}Let $\alpha$ be the composition
\[\wedge^i V_1 \to \wedge^i(V_0\otimes V_0) \to \wedge^i V_0 \otimes \Sym^i V_0,\]
and let $\beta$ be the composition
$$\wedge^i V_1 \otimes \Sym^{b-1-i} V_0 \to \wedge^i V_0 \otimes \Sym^{i} V_0 \otimes \Sym^{b-1-i}V_0 \to \wedge^i V_0 \otimes \Sym^{b-1} V_0$$
In the case of $V_1 = \wedge^2 V_0$, $\beta$ is $SL(V_0)$-equivariant and nonzero. Thus, to show \Cref{psidef}, it is sufficient to that the image of $\beta$ is contained in the kernel of the Koszul differential, $\ker(\delta)$ \eqref{Schurdef}. The map $\wedge^i(V_0\otimes V_0) \to \wedge^{i}V_0 \otimes \Sym^i V_0$ is given by multiplication from the differential graded algebra structure on the Koszul resolution interpreting $V_0 \otimes V_0$ as the elements in homological and internal degrees both being one (See \cite{Avramov}[Definition 1.3] for details on DG algebras). As $\wedge^2 V_0 \subset V_0 \otimes V_0$ are Koszul cycles and Koszul cycles are closed under multiplication in the DG algebra structure, we see that the image of $\alpha$ consists of cycles now in homological and internal degrees both $i$. As multiplication by elements of $S$ also preserves being a cycle, we see that the image of $\beta$ also consists of cycles now in homological degree $i$ and internal degree $b-1$. Thus, the image of $\beta$ is contained in $\ker{\delta}$, and so we obtain the map $\psi_{i,j}(V_1)$. Note further that by its construction, $\psi_{i,j}(V_1)$ is a restriction of $\psi_{i,j}(\wedge^2 V_0)$, that is that we have a commutative diagram: 
$$\begin{tikzcd}
\wedge^i V_1 \otimes \Sym^{b-1-i}V_0 \arrow[rd, "\psi_{i,j}(V_1)"] \arrow[rr, hook] & & \wedge^i(\wedge^2V_0) \otimes \Sym^{b-1-i}V_0 \arrow[ld, "\psi_{i,j}(\wedge^2 V_0)"]\\
& \mathbb{S}_{(i+j, 1^i)} V_0
\end{tikzcd}$$

\end{proof}

In the following proposition, we will use \Cref{psidef} to define a chain map 
\[\Sym^{b-1}(\varphi|_{V_1 \otimes S}) \to C_{\bullet}\]
for any choice of $V_1 \subset \wedge^2 V_0$.

\begin{proposition}\label{ischain}{Let $V_1 \subset \wedge^2 V_0$, and let \[f(V_1)_{i} = \psi_{i, b-1-i}(V_1): \wedge^iV_1 \otimes \Sym^{b-1-i}V_0 \to \mathbb{S}_{b-1, 1^i}V_0\] then $f(V_1)_{\bullet}: \Sym^{b-1}(\varphi|_{V_1) \to C_{\bullet}}$ is a chain map.}\end{proposition}
\begin{proof}
 To show that $f_{\bullet}(V_1)$ defines a map of chain complexes for any $V_1 \subset \wedge^2 V_0$, it is sufficient to show that $f(\wedge^2 V_0)_{\bullet}$ defines a map of chain complexes. To do this, we need to show that the diagram\\
$$\begin{tikzcd}\label{comm}
\wedge^i(\wedge^2V_0)\otimes \Sym^{b-1 - i} V_0 \arrow[r] \arrow[d, "f_i'"] & \wedge^{i-1}(\wedge^2 V_0) \otimes \Sym^{b-i} V_0 \otimes V_0\arrow[d, "f_{i-1}'\otimes \mbox{ id}"]\\
\mathbb{S}_{(b-1, 1^i)} V_0 \arrow[r] & \mathbb{S}_{(b-1, 1^{i-1})}\otimes V_0
\end{tikzcd}$$
commutes. All of the maps are $SL(V_0)$-equavariant. In particular, $f(\wedge^2 V_0)_i$ is just projection onto an irreducible factor. Further, $\mathbb{S}_{(b-1, 1^{i-1})} \otimes V_0$ has $SL(V_0)$-decomposition given by Pieri's rule:
$$\mathbb{S}_{(b-1, 1^{i-1})} \otimes V_0 \cong \mathbb{S}_{(b, 1^{i-1})}V_0 \oplus \mathbb{S}_{(b-1, 2, 1^{i-2})}V_0 \oplus \mathbb{S}_{(b-1, 1^i)}V_0$$

 Thus, if we can show that for the irreducible decomposition of $\wedge^i(\wedge^2 V_0) \otimes\Sym^{b-1-i}V_0$ into $SL(V_0)$ irreducible representations, $\mathbb{S}_{(b-1, 1^{i})}$ occurs with multiplicity one and the other two irreducible factors of  $\mathbb{S}_{(b-1, 1^{i-1})} \otimes V_0$ do not occur, then we have shown that the diagram commutes up to a scalar.

Let $e_0, ... , e_b$ be a basis for $V_0$. The highest weight of $\wedge{i}(\wedge^2V_0) \otimes \Sym^{b-1-i}V_0$ is $(b-1, 1^i)$ corresponding to the basis element $\bigwedge_{1\leq j\leq i}(e_0\wedge e_j) \otimes e_0^{b-1-i}$ (see \cite{FultonHarris}[Chapter 12] for an introduction to weight theory) . This is the only basis element of this weight, so $\wedge^i(\wedge^2V_0) \otimes \Sym^{b-1-i}V_0$ contains one copy of $\mathbb{S}_{(b-1, 1^i)}V_0$. $\wedge^i(\wedge^2V_0) \otimes \Sym^{b-1-i}V_0$ does not have any basis elements with weight $(b, 1^{i-1})$ or $(b-1, 2, 1^{i-2})$, so $\wedge^i(\wedge^2V_0) \otimes \Sym^{b-1-i}V_0$ does not have $\mathbb{S}_{(b, 1^{i-1})} V_0$ nor $\mathbb{S}_{(b-1, 2, 1^{i-2})} V_0$ as irreducible factors. 

We now show that these maps agree, by showing that they do so for $x = \bigwedge_{1\leq j\leq i}(e_0\wedge e_j) \otimes e_0^{b-1-j}$. The map $\mathbb{S}_{(b-1, 1^i)} V_0 \to \mathbb{S}_{(b-1, 1^{i-1})}V_0 \otimes V_0$ is given by restricting the comultiplication
$$\wedge^{i}V_0 \otimes \Sym^{b-1}V_0 \xrightarrow{d} \wedge^{i-1}V_0 \otimes V_0 \otimes \Sym^{b-1}V_0 = \wedge^{i -1 } V_0 \otimes \Sym^{b-1}V_0 \otimes V_0$$
Now, we compute
\[d(f(\wedge^2V_0)_i(x)) = d(e_1\wedge...\wedge e_i \otimes e_0^{b-1} + \sum_{0 \leq j \leq i} (-1)^j e_0\wedge ... \wedge \hat{e}_j \wedge ... \wedge e_i\otimes e_0^{b-2}\cdot e_j)\]
\[= d(\sum_{0 \leq j \leq i} (-1)^je_0 \wedge ... \wedge \hat{e}_j \wedge ... \wedge e_i \otimes e_0^{b-2}\cdot e_j)\]
\[= \sum_{0 \leq k < j \leq i} (-1)^{k + j} e_0\wedge...\wedge\hat{e}_k\wedge...\wedge \hat{e}_j\wedge...\wedge e_i \otimes (e_0^{b-2}\cdot e_j \otimes e_k - e_0^{b-2}\cdot e_k \otimes e_j)\]
Next,
\[f(\wedge^2V_0)_{i-1}(d(x))\]
\[ = f(\sum_{1 \leq j \leq i} (-1)^j(e_0 \wedge e_1)\wedge ... \wedge\widehat{(e_0 \wedge e_j)}\wedge ...\wedge(e_0\wedge e_i)\otimes (-e_0^{b-i}\otimes e_j + e_0^{b-1-i}\cdot e_j \otimes e_0))\]
\[ = \sum_{1 \leq j \leq i}(-1)^j e_1 \wedge .... \wedge \hat{e}_j \wedge... \wedge e_i \otimes e_0^i(-e_0^{b-i}\otimes e_j + e_0^{b-1-j} e_j \otimes e_0)\]
\[+\sum_{0 \leq k < j \leq i} [(-1)^{k + j} e_0\wedge...\wedge\hat{e}_k\wedge...\wedge \hat{e}_j\wedge...\wedge e_i\]
\[\otimes (e_0^{i-2}\cdot e_k (-e_0^{b-i}\otimes e_j + e_0^{b-1-i}\cdot e_j\otimes e_0) - e_0^{i-2}\cdot e_j(-e_0^{b-i}\otimes e_k + e_0^{b-1-i}\cdot e_k \otimes e_0))]\]
\[ = \sum_{1 \leq j \leq i}(-1)^j e_1 \wedge .... \wedge \hat{e}_j \wedge... \wedge e_i \otimes e_0^i(-e_0^{b-i}\otimes e_j + e_0^{b-1-j} e_j \otimes e_0)\]
\[+\sum_{0 \leq k < j \leq i} [(-1)^{k + j} e_0\wedge...\wedge\hat{e}_k\wedge...\wedge \hat{e}_j\wedge...\wedge e_i\]
\[\otimes(-e_0^{b-2}e_k\otimes e_j + e_0^{b-2}e_ke_j\otimes e_0 + e_0^{b-2}e_j\otimes e_k - e_0^{b-3}e_ke_j \otimes e_0)]\]
\[= \sum_{0 \leq k < j \leq i} (-1)^{k + j} e_0\wedge...\wedge\hat{e}_k\wedge...\wedge \hat{e}_j\wedge...\wedge e_i \otimes (e_0^{b-2}\cdot e_j \otimes e_k - e_0^{b-2}\cdot e_k \otimes e_j)
=d(f(\wedge^2V_0)_i(x))\]

\end{proof}

Next we will show that under suitable hypothesis on $V_1 \subset \wedge^2 V_0$ the chain map $f_{\bullet}(V_1)$ is injective by giving a general criteria for the maps $\psi_{i,j}(V_1)$ to be injective. To do so, we will first realize the maps $\psi_{i,j}$ cohomologically. Consider the exact sequence of sheaves on $\mathbb{P}(V_0)$.
\[\wedge^2 V_0 \otimes \Oo_{\Pp(V_0)} \to V_0 \otimes \Oo_{\Pp(V_0)}(1) \to \Oo_{\Pp(V_0)}(2) \to 0\]
For the rest of the paper fix the notation
\begin{equation}\label{tauto}K := \ker(V_0 \otimes \Oo_{\Pp(V_0)} \to \Oo_{\Pp(V_0)}(1)).\end{equation}
Note that $K = \Omega(1)$ where $\Omega$ is the cotangent sheaf on projective space. As the sequence of sheaves is exact, we have a surjection $\wedge^2 V_0 \otimes \mathbb{O}_{\mathbb{P}} \to K(1)$. If we take $\wedge^i$ of this map and then twist by $\mathbb{O}_{\mathbb{P}}(j)$, we obtain
\[\wedge^i(\wedge^2 V_0) \otimes \mathbb{O}_{\mathbb{P}}(j) \to \wedge^iK(i+j)\]
and on global sections this becomes
\[\wedge^i(\wedge^2 V_0) \otimes \Sym^j(V_0) \to \mathbb{S}_{i+j, 1^i}V_0\]
This map agrees with $\psi_{i,j}(\wedge^2 V_0)$ as it is $SL(V_0)$-equivariant and is nonzero. The map $\psi_{i,j}(V_1)$ is then global sections of the restriction 
\begin{equation}\label{cohomreform}\wedge^i(V_1) \otimes \Oo_{\Pp}(j) \to (\wedge^iK)(i+j) \end{equation}
 We will use the following to construct a complex of sheaves ending in $\wedge^{i}V_1 \otimes \Oo_{\Pp(V_0)} \to (\wedge^{i}K)(i)$ in order to use a hypercohomology spectral sequence to check the injectivity of $\psi_{i,j}(V_1)$.
\begin{lemma}\label{ses}
Suppose $\mathscr{F} \twoheadrightarrow \mathscr{G}$ is a surjective morphism of locally free sheaves on a Noetherian nonsingular quasiprojective variety $X$ over the field $k$, where $\mbox{rk } \mathscr{F} = f$ and $\mbox{rk } \mathscr{G} = g$. Then for any $1 \leq i \leq g$ there exists a long exact sequence of sheaves
$$0 \to \wedge^{f} \mathscr{F}\otimes \mathbb{S}_{(f-g, 1^{g-i})}\mathscr{G}^{\vee}\to ... \to \wedge^{g+2} \mathscr{F} \otimes \mathbb{S}_{(2, 1^{g-i})}\mathscr{G}^{\vee}$$
$$\to \wedge^{g+1} \mathscr{F} \otimes \mathbb{S}_{(1^{g+1-i})}\mathscr{G}^{\vee} \to \wedge^{i} \mathscr{F} \to \wedge^{i}\mathscr{G} \to 0$$
\end{lemma}
\begin{proof}
Consider $Y = Gr(i, \mathscr{G})$, the $i-{\mbox{th}}$ Grassmanian bundle of $\mathscr{G}$, which has projection $\pi: Y \to X$. We have the tautological sequence 
$$0 \to \mathscr{R} \to \pi^{*} \mathscr{G} \to \mathscr{Q} \to 0,$$
where $\mathscr{R}$ and $\mathscr{Q}$ are the tautological vector bundles of ranks $b+1-i$ and $i$-respectively. 
As $\mathscr{F} \to \mathscr{G}$ is surjective, we also have that the composition $\pi^{*}\mathscr{F} \to \pi^{*} \mathscr{G} \to \mathscr{Q}$ is surjective. As $\mbox{rk } \mathscr{Q} = i$, we have an Eagon-Northcott Resolution \cite{Positivity}[Theorem B.2.2 (The case EN0)] 
$$0 \to \wedge^{f} (\pi^* \mathscr{F}) \otimes (\Sym^{f-i}\mathscr{Q})^{\vee} \to ... \to \wedge^{i+1} (\pi^{*}\mathscr{F}) \otimes (\mathscr{Q})^{\vee} \to \wedge^{i} \mathscr{F} \to \wedge^{i} \mathscr{Q} \to 0$$
We are now going to get the desired resolution from taking the relative hypercohomology spectral sequence of this complex \cite{Weibel}[Proposition 5.7.9; Proposition 5.7.10]. The hypercohomology spectral sequence has $E_1$ page given by 
$$E_{1}^{p, q} = \begin{cases}(R^{q}\pi)(\wedge^{i} \mathscr{Q}) & p = -1\\ (R^{q}\pi)(\wedge^{i - p}(\pi^*\mathscr{F}) \otimes (\Sym^{-p}\mathscr{Q})^{\vee}) & -(f-i) \leq p \leq 0\\
0 & \mbox{otherwise}\end{cases}$$
As the complex is exact and bounded, we know that this spectral sequence abuts to $0$. We calculate the derived pushforwards using the projection formula \cite{Hartshorne}[Exercise III.8.3] and the Borel-Weil-Bott Theorem \cite{Weyman}[Theorems 4.1.4,4.1.9](Weyman 4.1.4 and 4.1.9). We have the following derived pushforward calculations
$$(R^q\pi)(\wedge^i \mathscr{Q}) = \begin{cases}
    \wedge^{i} \mathscr{G} & q = 0\\
    0 & \mbox{otherwise}
\end{cases}$$
For the other derived pushforward calculations, we can rewrite the vector bundle slightly
$$\wedge^{i-p} (\pi^*)\mathscr{F} \otimes (\Sym^{-p} \mathscr{Q})^{\vee} = \pi^*(\wedge^{i -p}\mathscr{F}) \otimes \mathbb{S}_{(0^{i-1}, p)}\mathscr{Q}$$
By the projection formula, we have that
\[(R^q\pi)(\pi^*(\wedge^{i-p}\mathscr{F}) \otimes \mathbb{S}_{(0^{i-1}, p)} \mathscr{Q}) = \wedge^{i-p}\mathscr{F} \otimes (R^q\pi)(\mathbb{S}_{(0^{i-1}, p)}\mathscr{Q})\]
To calculate $(R^q\pi)(\mathbb{S}_{(0^{i-1}, p)}\mathscr{Q})$, we need to apply Bott's algorithm to the tuple $(0^{i -1}, p, 0^{g-i})$. From this we see that
$$(R^{q}\pi)(\wedge^{i - p}(\pi^*\mathscr{F}) \otimes (\Sym^{-p}\mathscr{Q})^{\vee})$$
$$ = \begin{cases}
    \wedge^{i - p}\mathscr{F} \otimes \mathbb{S}_{(0^{i-1}, -1^{g-i}, p+g-i)} \mathscr{G} & q = g-i \mbox{ and}\\ 
     & g-i +1 \leq -p \leq f-i\\
    0 & \mbox{otherwise}
\end{cases}$$
Note that $\mathbb{S}_{(0^{i-1},-1^{g-i}, p + g -i)}\mathscr{G}) = (\mathbb{S}_{(-p-g+i, 1^{g-i})}\mathscr{G})^{\vee}$. Thus, the first page of the hypercohomology spectral sequence is \\
\adjustbox{scale = 0.97, center}{
\begin{tikzcd}
    \wedge^{f}\mathscr{F} \otimes (\Sym^{f-g}\mathscr{G})^{\vee} \arrow[r] & \dots \arrow[r]& \wedge^g \mathscr{F} \arrow[r] & 0 \arrow[r] & \dots \arrow[r] & 0  \arrow[r] & 0 \arrow[r] & 0\\
    0 \arrow[r] & \dots \arrow[r] & 0 \arrow[r] & 0 \arrow[r] & \dots \arrow[r] & 0  \arrow[r] & 0 \arrow[r] & 0\\
    \vdots\\
    0 \arrow[r] & \dots \arrow[r] & 0 \arrow[r] & 0 \arrow[r] & \dots \arrow[r] & 0  \arrow[r] & 0 \arrow[r] & 0\\
    0 \arrow[r] & \dots \arrow[r] & 0 \arrow[r] & 0 \arrow[r] & \dots \arrow[r] & 0  \arrow[r] & \wedge^i \mathscr{F} \arrow[r] & \wedge^i \mathscr{G}
\end{tikzcd} }
By the arrangement of 0's on the $E_1$ page, we see that on the $E_{g-i}$ page, there is a connecting homomorphism
\[\mbox{coker}(\wedge^{g+1} \mathscr{F} \otimes (\mathscr{G})^{\vee} \to \wedge^{g} \mathscr{F}) \to \ker(\wedge^{i} \mathscr{F} \to \wedge^{i}\mathscr{G})\]
Further, this is the only possible nonzero morphism after the $E_1$ page. Thus, we have a long exact sequence
$$0 \to \wedge^{f} \mathscr{F}\otimes \mathbb{S}_{(f-g, 1^{g-i})}\mathscr{G}^{\vee}\to ... \to \wedge^{g+2} \mathscr{F} \otimes \mathbb{S}_{(2, 1^{g-i})}\mathscr{G}^{\vee}$$
$$\to \wedge^{g+1} \mathscr{F} \otimes \mathbb{S}_{(1^{g+1-i})}\mathscr{G}^{\vee} \to \wedge^{i} \mathscr{F} \to \wedge^{i}\mathscr{G} \to 0$$
\end{proof}

We will now use \Cref{ses}, to give a criteria for the injectivity of the maps $\psi_{i,j}(V_1)$.
\begin{lemma}\label{psiinj}
Let $V_1 \subset \wedge^2 V_0$ such that $\mbox{dim}_{\mathbb{C}} V_1 = 2b-1$ and the middle homology of the sequence
\[V_1 \otimes S \to V_0 \otimes S \to S\]
is of finite length. Let $\psi_{i,j}(V_1)$ be as defined in \Cref{psidef}. If $i,j$ satisfy 
\[\begin{cases}
i+j < b & or\\
i = b, j= 0
\end{cases}\]
Then $\psi_{i,j}(V_1)$ is injective.
\end{lemma}
\begin{proof}
\noindent By hypothesis, we have a surjective map 
\[V_1 \otimes \Oo_{\Pp(V_0)} \to K(1),\]
where $K$ is the kernel in the tautological sequence \eqref{tauto}. If we apply \Cref{ses} to this surjection, we get a long exact sequence
\[0 \to \wedge^{2b-1}(V_1) \otimes \mathbb{S}_{(b-1, 1^{b-i})}((K(1))^{\vee}) \to ... \to \wedge^{b+1}(V_1) \otimes \mathbb{S}_{(1^{b+1 -i})} ((K(1))^{\vee})\] \[\to \wedge^i(V_1) \otimes \mathbb{O}_{\Pp(V_0)} \to \wedge^{i} (K(1))\]
Note that $\mathbb{S}_{(l, 1^{b-i})}((K(1)^{\vee})) = \mathbb{S}_{(0^{i-1}, -1^{b-i}, -l)}K(-(b+l-i))$
We now twist this sequence by $\mathbb{O}_{\Pp(V_0)}(j)$ to get the long exact sequence
\[0 \to \wedge^{2b-1}V_1 \otimes \mathbb{S}_{(0^{i-1}, -1^{b-i}, -(b-1))}K(-(2b-1-i-j)) \to ...\] \[\to \wedge^{b+l}V_1 \otimes \mathbb{S}_{(0^{i-1}, -1^{b-i}, -l)}K(-(b+l-i-j))
 \to ... \to \wedge^{b+1}V_1 \otimes \mathbb{S}_{(0^{i-1}, -1^{b+1-i})}(-(b+1-i-j))\] \[\to \wedge^i V_1 \otimes \mathbb{O}_{\Pp(V_0)}(b-1-i) \to \wedge^iK(b-1) \to 0\]
As this sequence is exact, the associated hypercohomology spectral sequence will again abut to 0. We wish to show the injectivity of the map
\[H^0(\wedge^i V_1 \otimes \mathbb{O}_{\Pp(V_0)}(b-1-i) \to \wedge^iK(b-1)),\]
so by the hypercohomology spectral sequence, it is enough to show the vanishing
\[H^{l-1}(\mathbb{S}_{(0^{i-1}, -1^{b-i}, -l)}K(-(b+l-i-j))) = 0\]
for $1 \leq l \leq b-1$. By the Borel-Weil-Bott theorem \cite{Weyman}[Theorems 4.1.4,4.1.9], we can determine for which $l$, $(\mathbb{S}_{(0^{i-1}, -1^{b-i}, -l)}K(-(b+l-i-j)))$ has nonzero cohomology groups. The only possibilities are the following three cases
$$\begin{cases}
H^{i-1} \neq 0 &  l = 2i+j-b \hspace{0.2cm} \& \hspace{0.2cm} i < b\\ 
H^{b-1} \neq 0 &  1 \leq i+j \leq l \\
H^{b} \neq 0 &  i+j <=1
\end{cases}$$
The latter two cases would give $H^{l-1} \neq 0$ if $l = b, b+1$ respectively, but this outside of the range $1 \leq l \leq b-1$ (Note this is why the assumption that $\dim V_1 = 2b-1$ is necessary in the hypothesis of the Lemma). In order for the first case to give $H^{l-1} \neq 0$, we would need that $l = i < b$. This would then give
\[i = 2i +j-b,\]
and so 
\[i + j = b.\]
By hypothesis, $i +j < b$ in this case, so the necessary vanishing occurs for all terms and we are done.
\end{proof}
\begin{corollary}\label{finj}
    Let $V_1 \subset \wedge^2 V_0$ satisfy the hypotheses of \Cref{psiinj}. Then $f(V_1)_{\bullet}$ is injective.
\end{corollary}
\begin{proof}
    $f(V_1)_i = \psi_{i,b-1-i}(V_1)$ by definition, so each map is injective making $f(V_1)_{\bullet}$ injective.
\end{proof}
\begin{remark} By the proof of \Cref{psiinj}, 
\[\coker(f(V_1)_i) = \wedge^{b+i+1}V_1 \otimes H^{i-1}(\Ss_{(0^{i-1}, -1^{b-i}, -(i-1))}K(-i))\] \[= \wedge^{b+i +1 V_1} \otimes \Ss_{(-1^b, -(i-1))}\] 
\[=(\wedge^{b-i-2}V_1)^* \otimes (\Sym^{i}V_0)^*\otimes \mbox{det}(V_0)^*\]
which is a term of the complex $\wedge^{b-2} \varphi^*$. This suggests that the resolution of $m^{b-1}$ can be realized as a mapping cone of $\Sym^{b-1} \varphi$ and $\wedge^{b-2} \varphi^*$.\end{remark}

\section{Bernstein-Gel'fand-Gel'fand Correspondence and Modules Generated in Lowest Degree}\label{BGG}
We will now use the Bernstein-Gel'fand-Gel'fand correspondence \cite{syzygies}[Chapter] (henceforth referred to as the BGG correspondence), which is an equivalence of categories between the category of linear complexes over the symmetric algebra, $\Sym(V_0)$, and the category of graded modules over the exterior algebra, $E :=\wedge(V_0^*)$. Here we give $E$ the grading given by $E_{-j} = \wedge^j(V_0^*)$. For $V_1 \subset \wedge^2 V-0$, let $P(V_1)$ be the module over $E$, that corresponds to the complex $\Sym^{b-1}(\varphi|_{V_1 \otimes S})$ and let $\hat{P}(V_1)$ be the module that corresponds to $\wedge^{b-1}(\varphi|_{V_1 \otimes S}^*)$. Via the BGG correspondence, we have that $P$ can be decomposed into its graded pieces as 
$$P(V_1) = \oplus_{0 \leq j \leq b-1} P_j = \oplus_{0 \leq j \leq b-1} \wedge^{j}V_1 \otimes \mbox{Sym}^{b-1-j}V_0.$$ 
Likewise $\hat{P}(V_1)$ has graded decomposition given by
$$\hat{P}(V_1) = \oplus_{-b-1 \leq j \leq 0} \wedge^{-j}V_1^* \otimes (\mbox{Sym}^{b-1+j}V_0)^*$$
Let $Q$ be the module over $E$ that corresponds to $C_{\bullet}$, 
the resolution of $m^{b-1}$
introduced in \Cref{embedding}. 
Corresponding to the chain map, $f(V_1)_{\bullet}: \mbox{Sym}^{b-1}(\varphi|_{V_1 \otimes S}) \to C_{\bullet}$, introduced in \Cref{embedding}, we have with an abuse of notation a graded map of $E$-modules $f(V_1):P(V_1) \to Q$. Under the hypotheses of \Cref{finj}, the map of graded modules $f(V_1)$ is also injective. Via the BGG correspondence, the conclusion of \Cref{main1} can be rephrased as 
$$P \cong \hat{P}(-(b-1))$$
This is the version that will be proven in \Cref{proofsec}. The main results of this section are that both $P(V_1)$ and $\hat{P}(V_1)$ are generated in a single degree under a set of assumptions $V_1 \subset \wedge^2 V_0$. 
\begin{corollary}
Let $V_1 \subset \wedge^2 V_0$ such that $\mbox{dim}_{\mathbb{C}} V_1 = 2b-1$ and the middle homology of the sequence
\[V_1 \otimes S \to V_0 \otimes S \to S\]
is of finite length. Then $\hat{P}(V_1)$ is generated in degree $0$, that is by $(\mbox{Sym}^{b-1}V_0)^*$.
\end{corollary}
\begin{proof}
By \Cref{finj}, $f(V_1)$ is injective. For $i = 0$, $f(V_1)_0$ is the identity on $\Sym^{b-1} V_0$. Further, 
\[\mbox{dim}\quad (V_1 \otimes \Sym^{b-2}V_0) = (2b-1) * \binom{2b-2}{b} = \mbox{dim} (\Ss(b-1, 1)V_0)\]
by the hook length formula \cite{FultonHarris}[Exercise A.30]. As $f_1$ is injective, it is also an isomorphism. As $C_{\bullet}$ is exact and $f(V_1)_{\bullet}$ is injective, $f(V_1)_{\bullet}$ is an embedding of $\mbox{Sym}^{b-1}(\varphi|_{V_1})$ into the linear strand of the resolution of $H_0(\mbox{Sym}^{b-1}\varphi|_{V_1})$. By \cite{syzygies}[Theorem 7.4 part 2], $\hat{P}$ is generated in degree $0$.
\end{proof}
We will now work towards criteria for $P(V_1)$ to be generated in lowest degree in several steps. In each degree $1 \leq j \leq b-1$, we need to show that the multiplication map
$$m_{-j}: \wedge^j V_1 \otimes \mbox{Sym}^{b-1-j} V_0 \otimes V_0^* \to \wedge^{j-1}V_1 \otimes \mbox{Sym}^{b-j}V_0$$
is surjective. Note that this map can be factored as
\[\wedge^jV_1 \otimes \mbox{Sym}^{b-1-j}V_0\otimes V_0^* \xrightarrow[\Delta \otimes\mbox{id}]{} \wedge^{j-1}V_1 \otimes V_1 \otimes V_0^* \otimes \mbox{Sym}^{b-1-j}(V_0)\]
\[\xrightarrow[\mbox{id}\otimes \theta \otimes \mbox{id}]{} \wedge^{j-1}V_1 \otimes V_0 \otimes \mbox{Sym}^{b-1-j} \xrightarrow[\mbox{id}\otimes m]{} \wedge^{j-1}V_1 \otimes \mbox{Sym}^{b-j}V_0\]
where $\Delta$ is comultiplication on the exterior algebra, $m$ is multiplication on the symmetric algebra, and $\theta$ is the composition
\[ V_1 \otimes V_0^* \xrightarrow[\varphi|_{V_1 \otimes S} \otimes \id]{} V_0 \otimes V_0 \otimes V_0^* \to V_0\]
As the multiplication map 
\[\mbox{Sym}^{b-1-j}V_0 \otimes V_0 \to \mbox{Sym}^{b-j}V_0\]
is surjective. It is sufficient to show the stronger statement that
$$\wedge^j V_1 \otimes V_0^* \to \wedge^{j-1} V_1 \otimes V_0$$
is surjective for $1 \leq j \leq b-1$. We will first prove this surjectivity in the case of $V_1 = \wedge^2 V_0$ before specializing to other choices of $V_1 \subset \wedge^2 V_0$.
\begin{lemma}\label{gen0}
 For each $1 \leq j \leq b-1$, the composition
$$\wedge^j (\wedge^2V_0) \otimes  V_0^* \xrightarrow[\Delta \otimes \mbox{id}]{} \wedge^{j-1}(\wedge^2 V_0) \otimes \wedge^2V_0 \otimes V_0^*$$
$$\xrightarrow[\mbox{id} \otimes \theta]{} \wedge^{j-1}(\wedge^2 V_0) \otimes V_0$$
is surjective.
\end{lemma}
\begin{proof}
First note that in this case,  $\theta$ is an $SL(V_0)$-equivariant projection onto an irreducible factor:
\[\theta: \wedge^2 V_0 \otimes V_0^* \to V_0\]
The maps
$$\wedge^j (\wedge^2 V_0) \otimes V_0^* \to \wedge^{j-1}V_0 \otimes V_0$$
 can be interpreted as maps between cohomology groups of sheaves on projective space as follows. As before, let $K$ be the kernel in the tautological sequence
\[0 \to K \to V_0 \otimes \Oo_{\Pp(V_0)} \to \Oo_{\Pp(V_0)}(1)\]
on $\mathbb{P}(V_0)$. The sheaves $\wedge^i K$ are twists of syzygy sheaves for the Koszul complex on $\Pp(V_0)$, so in particular we have a short exact sequence
\begin{equation}\label{ses2} 0 \to \wedge^2 K \to \wedge^2 V_0 \otimes \Oo_{\Pp(V_0)} \to K(1) \to 0\end{equation}
This yields a long exact sequence
\[0 \to \wedge^j(\wedge^2K) \to \Sym^j(\wedge^2V_0 \otimes \Oo_{\Pp(V_0)} \to K(1)) \to 0\]
that is a long exact sequnce
\[0 \to \wedge^j(\wedge^2 K) \to \wedge^j(\wedge^2 V_0) \otimes \Oo_{\Pp(V_0)} \to \wedge^{j -1}(\wedge^2 V_0) \otimes K(1) \to \cdots \to \Sym^jK(j) \to 0\]
Tensoring this sequence by $\Oo_{\Pp(V_0)}(-b-2)$ gives a long exact sequence
\[0 \to \wedge^j\wedge^2K(-b-2) \to \wedge^j(\wedge^2 V_0) \otimes \Oo_{\Pp(V_0)}(-b-2) \to \wedge^{j -1}(\wedge^2 V_0) \otimes K(-b-1)\]
\[ \to \cdots \Sym^jK(-b-2+j) \to 0\]
We now show that, $$H^b(\wedge^j(\wedge^2 V_0) \otimes \Oo_{\Pp(V_0)}(-b-2) \to \wedge^{j -1}(\wedge^2V_0) \otimes K(-b-1))$$
is the desired map
\[\wedge^j (\wedge^2 V_0) \otimes V_0^* \to \wedge^{j-1}V_0 \otimes V_0\]
(tensored with $\wedge^{b+1}(V_0)^* = (\mbox{det}V_0)^*$ and up to a scalar).
First we have that, 
$$H^b(\Oo_{\Pp(V_0)}(-b-2)) = V_0^* \otimes \mbox{det}(V_0)^*$$
To calculate the cohomology of $K(-b-1)$, we need to apply the Borel-Weil-Bott theorem \cite{Weibel}[Theorems 4.1.4,4.1.9] to the tuple $(-b-1, 1, 0^{b-1})$, which results in the tuple $(0, -1^b)$ after $b$ transpositions, so 
$$H^b(K(-b-1)) = \mbox{det}(V_0)^* \otimes V_0$$
Further, the map
$$\wedge^j(\wedge^2 V_0) \otimes \Oo_{\Pp(V_0)}(-b-2) \to \wedge^{j -1}(\wedge^2V_0) \otimes K(-b-1)$$
factors as 
\[\wedge^j(\wedge^2 V_0) \otimes \Oo_{\Pp(V_0)}(-b-2) \xrightarrow[\Delta \otimes \mbox{id}]{} \wedge^{j-1}(\wedge^2V_0)\otimes \wedge^2V_0 \otimes \Oo_{\Pp(V_0)}(-b-2)\]
\[\xrightarrow[\mbox{id} \otimes \theta']{} \wedge^{j -1}(\wedge^2V_0) \otimes K(-b-1)\]
Here $\theta': \wedge^2 V_0 \otimes \Oo_{\Pp(V_0)}(-b-2) \to K(-b-1)$ is the map from \eqref{ses2} twisted by $\Oo_{\Pp(V_0)}(-b-2)$. Applying the long exact sequence on cohomology to \eqref{ses2}, $H^b(\theta') \neq 0$. Hence $H^b(\theta')$ is a nonzero $SL(V_0)$-equavariant map
\[\wedge^2V_0 \otimes V_0^* \to V_0,\]
so it must be projection onto an irreducible factor as a representation of $SL(V_0)$, i.e. agrees with $\theta$ up to a nonzero scalar. Thus,
 $$H^b(\wedge^j(\wedge^2 V_0) \otimes \Oo_{\Pp(V_0)}(-b-2) \to \wedge^{j -1}(\wedge^2V_0) \otimes K(-b-1))$$
is the desired map. We now show the desired surjectivity by using the hypercohomology sequence associated to the resolution of sheaves above. As $H^{b+1}(\wedge^j(\wedge^2 K)(-b-2)) = 0$ ($\Pp(V_0)$ has dimension $b$), it is enough to show the vanishing
\[H^{b-l+2}(\wedge^{j-l}(\wedge^2V_0) \otimes \mbox{Sym}^lK(l-b-2))\]
or equivalently
\[H^{b-l+2}( \mbox{Sym}^lK(l-b-2))\]
for $2 \leq l \leq j$. Note that as $j \leq b-1$, we have that $b-l+2 \geq 3$. We now compute the possible nonzero cohomology groups by applying the Bott algorithm \cite{Weyman}[Theorems 4.1.4,4.1.9] to the tuple (-(b+2-l), l, 0). There are two possible cases
\[\begin{cases}
    H^1(\Sym^lK(l-b-2))\neq 0 & l \geq b-1\\
    H^b(\Sym^l K(l-b-2)) \neq 0 & l \leq 1
\end{cases}\]

As $b-l+2 \geq 3$, the first case cannot give $H^{b-l+2} \neq 0$. The latter case is outside of the range for $l$, so we have the necessary vanishing. 
\end{proof}

We now show that certain subspaces of $\wedge^2V_0$ also satisfy the conclusion of \Cref{gen0}.
\begin{lemma}\label{gen1}
Let $V_1 \subset \wedge^2V_0$ such that if we restrict the map 
$$\wedge^2V_0 \otimes V_0^* \to V_0$$
to $V_1 \otimes V_0^*$ it remains surjective. Then the restriction of the map from \Cref{gen0}
$$ \wedge^{j}(V_1) \otimes V_0^* \to \wedge^{j-1}(\wedge^2 V_1) \otimes V_0$$
is also surjective.
\end{lemma}

\begin{proof}
Consider the diagram
\[\begin{tikzcd}
\wedge^j V_1 \otimes V_0^* \arrow[hookrightarrow]{rr} \arrow[d, "\alpha_1'"] & & \wedge^j (\wedge^2V_0) \otimes V_0^* \arrow[d, "\alpha_1"]\\
\wedge^{j-1} V_1 \otimes V_1 \otimes V_0^* \arrow[hookrightarrow]{r} \arrow[d, "\alpha_2'"]& \wedge^{j-1}V_1 \otimes \wedge^2 V_0 \otimes V_0^* \arrow[hookrightarrow]{r} \arrow[dl, "\alpha_2''"] &\wedge^{j-1}(\wedge^2 V_0) \otimes \wedge^2V_0 \otimes V_0^* \arrow[d, "\alpha_2"]\\
\wedge^{j-1}V_1 \otimes V_0 \arrow[hookrightarrow]{rr} & & \wedge^{j-1}(\wedge^2V_0) \otimes V_0 
\end{tikzcd}\]

In this diagram, all of the horizontal maps are inclusions, $\alpha_1'$ is a restriction of $\alpha_1$ under the inclusion, and $\alpha_2'$ and $\alpha_2''$ are restrictions of $\alpha_2$ under the inclusions. By \Cref{gen0}, $\alpha_2 \circ \alpha_1$ is surjective. We will show this surjectivity extends to $\alpha_2' \circ \alpha_1'$.

As a property of comultiplication on exterior algebras 
\[\alpha_1^{-1}(\wedge^{j-1}V_1 \otimes \wedge^2 V_0 \otimes V_0^*) = \wedge^jV_1 \otimes V_0^*\]
Further, 
\[\alpha_2^{-1}(\wedge^{j-1}V_1 \otimes V_0) = \wedge^{j-1}V_1 \otimes \wedge^2V_0 \otimes V_0^*\]
As $V_1 \otimes V_0^* \to V_0$ is surjective,
\[\alpha_2'(\wedge^{j-1}V_1 \otimes V_1 \otimes V_0) = \alpha_2''(\wedge^{j-1}V_1 \otimes \wedge^2V_0 \otimes V_0^*) = \wedge^{j-1}V_1 \otimes V_0\]
Putting this together with the surjectivity of $\alpha_2\circ \alpha_1$
\[\wedge^{j-1}V_1 \otimes V_0 = \wedge^{j -1}V_1 \otimes V_0 \cap \alpha_2\circ \alpha_1(\wedge^j(\wedge^2V_0)\otimes V_0^*)\]
\[= \alpha_2(\wedge^{j-1}V_1 \otimes \wedge^2V_0\otimes V_0^* \cap \alpha_1(\wedge^j(\wedge^2V_0)\otimes V_0^*))\]
\[= \alpha_2(\wedge^{j-1}V_1 \otimes V_1 \otimes V_0^* \cap \alpha_1(\wedge^j(\wedge^2 V_0)\otimes V_0^*))\]
\[= \alpha_2(\alpha_1(\wedge^j V_1 \otimes V_0^*))\]
\end{proof}

We now show that the hypotheses $V_1\subset \wedge^2 V_0$ repeated throughout the paper satisfy the hypotheses of \Cref{gen1}.
\begin{lemma}
Let $V_1 \subset \wedge^2 V_0$ such that the sequence of sheaves on $\mathbb{P}(V_0)$
\[V_1 \otimes \Oo_{\Pp(V_0)} \to V_0 \otimes \Oo_{\Pp(V_0)}(1) \to \Oo_{\Pp(V_0)}(2) \to 0\]
is exact. Then the composition 
\[V_1 \otimes V_0^* \to V_0 \otimes V_0 \otimes V_0^* \to V_0\] 
is surjective.
\end{lemma}
\begin{proof}
    First consider the composition, 
    \[V_1 \otimes V_0^* \otimes \mathbb{O}_{\mathbb{P}(V_0)} \to V_0 \otimes V_0^* \otimes \mathbb{O}_{\mathbb{P}(1)} \to \mathbb{O}_{\mathbb{P}(V_0)}(1) \]
    The composition of this map of sheaves is surjective if and only if the cokernel of the corresponding map of free modules has finite length. As the corresponding map of free modules is linear, this is equivalent to the composition of global sections being surjective. On global sections, the composition is
    \[V_1 \otimes V_0^* \to V_0\otimes V_0 \otimes V_0^* \to V_0,\]
    which is what we are trying to show is surjective.
    By assumption, the map 
    \[V_1 \otimes V_0^* \otimes \mathbb{O}_{\mathbb{P}(V_0)}  \to \mathbb{O}_{\mathbb{P}(V_0)}(1)\] factors as a sequence of maps
    \[V_1 \otimes V_0^* \otimes \mathbb{O}_{\mathbb{P}(V_0)} \twoheadrightarrow V_0^* \otimes K(1) \to \mathbb{O}_{\mathbb{P}(V_0)}(1)\]
    where $K$ is again the kernel in the tautological sequence \eqref{tauto}. Thus, it is sufficient to show that the map
    \[V_0^*\otimes K(1) \to \mathbb{O}_{{\mathbb{P}(V_0)}}(1)\]
    is surjective. If we tensor this map by $\mathbb{O}_{\mathbb{P}(V_0)}(e)$ for $e \geq 0$ and then take global sections, the maps are $SL(V_0)$-equivariant projections onto an irreducible factor
    \[V_0^* \otimes \mathbb{S}_{(e+1, 1)}V_0 \to \mbox{Sym}^{e+1}V_0\]
    each of which is surjective, so therefore the map of sheaves is also surjective.
\end{proof}

\section{Morphisms of $E$-modules via Hermite Reciprocity}\label{proofsec}
In this section, we start by defining a chain map $\wedge^{b-1}((\varphi|_{V_1 \otimes S})^*) \to C_{\bullet}$ for the case of $\dim V_1 = 2b-1$ in order to get a map of $E$-modules $\hat{P}(-b-1) \to Q$. When we showed that $f_{\bullet}(V_1)$ is a chain map, it was convenient to use the fact that $\Sym^{b-1}(\varphi|_{V_1 \otimes S})$ is a subcomplex of $\Sym^{b-1}(\varphi)$. Now, $\wedge^{b-1}((\varphi|_{V_1 \otimes S})^*)$ is instead a quotient of $\wedge^{b-1}(\varphi^*)$. We can get around this issue with the following lemma.

\begin{lemma}\label{wedgereform}
Let $V_1 \subset \wedge^2 V_0$ be of dimension $2b-1$, then $\wedge^{b-1}((\varphi|_{V_1 \otimes S})^*)$ is a subcomplex of $\wedge^{\binom{b}{2}}(\varphi^*)$.
\end{lemma}
\begin{proof}
    We use the identification \Cref{wedgedual}, to identify the complex $\wedge^{b-1}((\varphi|_{V_1 \otimes S})^*)$ with 
    \[\wedge^{2b-1}V_1 \otimes (\Sym^{b-1} V_0)^*\otimes S \to \wedge^{2b-2} V_1 \otimes (\Sym^{b-2} V_0)^* \otimes S \to \cdots\]
    \[\to \wedge^{b+1} V_1 \otimes V_0^*\otimes S  \to \wedge^b V_1 \otimes S\]
    and identify the complex $\wedge^{\binom{b}{2}}(\varphi)$ with
    \[\wedge^{\binom{b+1}{2}}(\wedge^2 V_0) \otimes (\Sym^{\binom{b}{2}} V_0)^*\otimes S \to \wedge^{\binom{b+1}{2}-1} (\wedge^2 V_0) \otimes (\Sym^{\binom{b}{2}-1} V_0)^* \otimes S \to \cdots\]
    \[\to \wedge^{b+1} (\wedge^2 V_0) \otimes V_0^*\otimes S  \to \wedge^b (\wedge^2 V_0) \otimes S\]
    The inclusion of complexes arises from the inclusion 
    \[\wedge^{b+i}(V_1) \otimes \Sym^iV_0^* \subset \wedge^{b+i}(\wedge^2V_0)\otimes \Sym^i V_0^*\]
    Additionally, the differentials of the complex $\wedge^{b-1}((\varphi|_{V_1 \otimes S})^*)$ arise as restrictions of the differentials of the complex $\wedge^{\binom{b}{2}}(\varphi^*)$, i.e. we have a commutative diagram
    \[\begin{tikzcd}
    \wedge^{b+i} V_1 \otimes \Sym^i V_0^*\otimes S \arrow[r] \arrow[d] & \wedge^{b+i} (\wedge^2(V_0)) \otimes \Sym^i V_0^*\otimes S \arrow[d]\\
    \wedge^{b+i-1}V_1 \otimes V_1 \otimes V_0^* \otimes \Sym^{i-1}V_0^* \otimes S \arrow[r] \arrow[d]  & \wedge^{b+i-1}(\wedge^2V_0) \otimes \wedge^2V_0 \otimes V_0^* \otimes \Sym^{i-1}V_0^* \otimes S \arrow[d]\\
    \wedge^{b+i-1}V_1 \otimes \Sym^{i-1}V_0^* \arrow[r] & \wedge^{b+i-1}(\wedge^2 V_0) \otimes \Sym^{i-1}V_0^*
    \end{tikzcd}\]
Here the horizontal arrows are inclusions, and the composition of vertical maps are the differentials for the complexes as given by \eqref{wedgedualdiff}
\end{proof}

Let $T_{\bullet}$ be the truncation of $\wedge^{\binom{b}{2}}(\varphi^*)$ consisting of the terms
\[\wedge^{b+j}(\wedge^2 V_0)\otimes \Sym^jV_0^*\]
By the proof of \Cref{wedgereform}, we have that $\wedge^{b-1}(\varphi|_{V_1}^*)$ is a subcomplex of $T_{\bullet}$ for $V_1$ of dimension $2b-1$. We now define a map of chain complexes, $g: T_{\bullet} \to C_{\bullet}$, which preserves the internal degree of $S$ but shifts the homological degree by $b-1$. We will later restrict $g_{\bullet}$ to $\wedge^{b-1}(\varphi|_{V_1}^*)$ for $V_1$ of dimension $2b-1$. In order for $g_{\bullet}$ to be more cleanly described, we will tensor the terms of $T_{\bullet}$ by $\mbox{det} V_0^*$ (this will also make the maps $GL(V_0)$-equavariant).  Let $g_i$ be the composition
\[\wedge^{b+j}(\wedge^2 V_0) \otimes S^jV_0^* \otimes \det V_0^* \to \wedge^{b+j}(\wedge^2 V_0)\otimes S^{b-1-j}V_0 \otimes S^{b-1-j}V_0^* \otimes S^jV_0^* \otimes \det V_0^*\]
\[\to \wedge^{b+j} (\wedge^2 V_0) \otimes \Sym^{b-1-j}V_0 \otimes \mathbb{S}_{b,1^b}V_0^*\]  \[\xrightarrow[]{\id \otimes \psi(\wedge^2V_0)_{b,0}^*} \wedge^{b+j}(\wedge^2V_0) \otimes \Sym^{b-1-j}V_0 \otimes \wedge^b(\wedge^2V_0)^*\]
\[\to \wedge^{j}(\wedge^2 V_0) \otimes \Sym^{b-1-j}V_0 \xrightarrow[]{f(\wedge^2 V_0)_j} \mathbb{S}(b-1, 1^j)V_0\]
\begin{lemma}
$g_{\bullet}$ is a chain map
\end{lemma}
We need to show the commutativity of the diagram
\[\begin{tikzcd}
    \wedge^{b+j}(\wedge^2 V_0) \otimes S^j(V_0)^* \otimes \det(V_0)^* \arrow[d] \arrow[r] & \wedge^{b+j-1}(\wedge^2 V_0) \otimes S^{j-1}(V_0)^* \otimes \det(V_0)^* \otimes V_0 \arrow[d]\\
    \mathbb{S}_{(b-1, 1^j)} V_0 \arrow[r] & \mathbb{S}_{(b-1, 1^{j-1})} V_0 \otimes V_0    
\end{tikzcd}\]
\begin{proof}
If we tensor this diagram by $\mbox{det}(V_0)$ (which does not affect whether or not the diagram is commutative), we obtain
\[\begin{tikzcd}
    \wedge^{b+j}(\wedge^2 V_0) \otimes S^j(V_0)^* \arrow[d] \arrow[r] & \wedge^{b+j-1}(\wedge^2 V_0) \otimes S^{j-1}(V_0)^* \otimes V_0 \arrow[d]\\
    \mathbb{S}_{(b, 2^j, 1^{b-j})} V_0 \arrow[r] & \mathbb{S}_{(b, 2^{j-1}, 1^{b-j+1})} V_0 \otimes V_0    
\end{tikzcd}\]
Checking the commutativity of this diagram will be very similar to the proof of \Cref{ischain}. All of the maps are also $SL(V_0)$-equivariant, so the first step is to check that among the $SL(V_0)$ irreducible factors of $\mathbb{S}_{(b, 2^{j-1}, 1^{b-j+1})} V_0 \otimes V_0$ only $\mathbb{S}_{b, 2^j, 1^{b-j}} V_0$ appears (and with a multiplicity of 1) in the irreducible decomposition for $\wedge^{b+j}(\wedge^2 V_0) \otimes S^j(V_0)^* \otimes \det(V_0)^*$. $\mathbb{S}_{(b, 2^{j-1}, 1^{b-j+1})} V_0 \otimes V_0$ decomposes as a direct sum 
\[\mathbb{S}_{(b, 2^{j}, 1^{b-j})} V_0 \oplus \mathbb{S}_{(b+1, 2^{j-1}, 1^{b-j+1})}V_0 \oplus \mathbb{S}_{(b, 3, 2^{j-2}, 1^{b-j})}V_0\]
$\wedge^{b+j}(\wedge^2 V_0) \otimes S^j(V_0)^*$ does not contain any vectors of weight $(b+1, 2^{j-1}, 1^{b-j+1})$, so\\ $\mathbb{S}_{(b+1, 2^{j-1}, 1^{b-j+1})}V_0$ is not part of the irreducible decomposition. To show that $\mathbb{S}_{(b, 2^{j}, 1^{b-j})} V_0$ appears with multiplicity one and $\mathbb{S}_{(b, 3, 2^{j-2}, 1^{b-j+1})}V_0$ does not appear, requires a more delicate plethysm argument. A description of the irreducible factors of $\wedge^{b+j}(\wedge^2 V_0)$ is given by \cite{Macdonald}[Chapter 1, Section 8, Exercise 6c]. The irreducible factors correspond to Young diagrams with $2(b+j)$ boxes, with at most $b+1$ rows, and such that in each hook of the diagram (a hook is formed by taking all of the boxes to the right or below a box on the diagonal), there is one more box in the column of the hook than in the row of the hook. Further, each such irreducible factor occurs with multiplicity one. By Pieri's rule, the irreducible factors for $\wedge^{b+j}(\wedge^2 V_0) \otimes S^j(V_0)^*$, correspond to the Young diagrams corresponding to the irreducible factors for $\wedge^{b+j}(\wedge^2 V_0)$ with $j$-boxes removed, no two removed from the same column, and such that the result is still a Young diagram.

$\mathbb{S}_{(b, 2^{j}, 1^{b-j})} V_0$ corresponds to a Young diagram with the outermost hook having a row of $b$ boxes and a column of $b+1$ boxes, and the inner hook having a row of $1$ box and a column of $j$-boxes. The only way to obtain this Young diagram by removing $j$ boxes from $j$ different columns of a Young diagram corresponding to a factor for $\wedge^{b+j}(\wedge^2 V_0)$ is to have a box from each of the second through $j+1$-st columns of the Young diagram for the partition $(b, j+1, 2^{j}, 1^{b-j+1})$. Thus, $\mathbb{S}_{(b,2^j, 1^{b-j})}V_0$ occurs with multiplicity one for the irreducible decomposition of $\wedge^{b+j}(\wedge^2 V_0) \otimes S^j(V_0)^*$.

$\mathbb{S}_{(b, 3, 2^{j-2}, 1^{b-j+1})}V_0$ corresponds to a Young diagram where the outermost hook has a row of $b$ boxes and a column of $b+1$ boxes, and the inner row has a row of $2$ boxes and a column of $j$ boxes. To attempt to add $j$-boxes so that the hooks are of the correct form, we need to add a box to each of the fourth to $j-1$-st columns, so that now the inner hook is of the correct form with a row of $j-1$ boxes and a column of $j$ boxes. We still need to add $4$ more boxes. We can't add any to the outermost hook, as the first column already has $b+1$ boxes and we can't have more than $b+1$-rows. We can add $2$ boxes to the second hook by adding a box to the second column and the $j$-th column keeping the hook in the correct form. Adding additional boxes to this hook would require adding additional boxes to the second column, which is not allowed, so we still have two additional boxes to add. Adding boxes to a new third hook, would require adding two boxes to the third column, which is also not allowed. Thus, $\mathbb{S}_{(b, 3, 2^{j-2}, 1^{b-j+1})}V_0$ does not occur in the irreducible decomposition for $\wedge^{b+j}(\wedge^2 V_0) \otimes S^j(V_0)^*$.

We conclude that the diagram commutes up to a scalar. To finish showing $g_{\bullet}$ is a chain map, we need to check that it commutes for a single element. This is a very tedious computation similar to the one done at the end of the proof of \Cref{ischain}, so we omit it here.
\end{proof}

We are now ready to give the proof of \Cref{main1}.
\begin{proof}[Proof of \Cref{main1}]
To show that 
\[\Sym^{b-1}(\varphi) \cong \wedge^{b-1}(\varphi^*),\] 
it is enough to show that the corresponding graded modules over $E$ are isomorphic, i.e.
\[P \cong \hat{P}(b-1)\]
As noted before, the injective map of chain complexes
\[f(V_1)_{\bullet}: \Sym^{b-1}(\varphi) \to C_{\bullet}\]
corresponds to an injective map of graded $E$-modules $f(V_1): P \to Q$. If we restrict the map 
\[g: T_{\bullet} \to C_\bullet\]
then we get a chain map
\[(g_{\bullet})|_{\wedge^{b-1}(\varphi^*)}: \wedge^{b-1}(\varphi^*) \to C_\bullet \]
and a corresponding map of graded $E$-modules $g: \hat{P}(b-1) \to C_{\bullet}$
As the Hilbert functions for $P$ and $\hat{P}(b-1)$ are the same and $f(V_1)$ is injective, it is enough to show that
\[g(\hat{P}(b-1)) \supseteq f(V_1)(P)\]
As $P$ is generated in a single degree, we only need to show this containment in that degree, that is
\[g_{b-1}(\wedge^{2b-1}V_1 \otimes \Sym^{b-1}V_0^*) \supseteq f(V_1)_{b-1}(\wedge^{b-1} V_1)\]
To accomplish this, we are going to show that if we restrict $g_{b-1}$ to $\wedge^{2b-1}V_1 \otimes \Sym^{b-1}V_0^*$ then it can be factored as
\begin{equation}\label{containmentfactor}\begin{tikzcd}
\wedge^{2b-1}V_1 \otimes \Sym^{b-1}V_0^* \arrow[d, "h"] \arrow[dr, "g_{b-1}"]\\
\wedge^{b-1} V_1 \arrow[r, "f(V_1)_{b-1}"] & \mathbb{S}_{b-1, 1^{b-1}}
\end{tikzcd}\end{equation}
where $h$ is given by the composition
\[\wedge^{2b-1}V_1 \otimes \Sym^{b-1}V_0^* \xrightarrow[]{\id \otimes \psi(V_1)_{b,0}^*} \wedge^{2b-1}V_1 \otimes \wedge^b V_1^* \to \wedge^{b-1}V_1\]
By \Cref{psiinj}, $\psi_{b,0}$ is injective. As 
\[\mbox{dim } \Sym^{b-1} V_0 \cong \wedge^{b-1} V_1,\]
$\psi_{b,0}$ is an isomorphism. Thus, $h$ is also an isomorphism, so we are done if we can show the commutativity of \eqref{containmentfactor}. To prove this commutativity, consider the diagram 
\[\begin{tikzcd}
\wedge^{2b-1}(\wedge^2 V_0) \otimes \Sym^{b-1}(V_0)^* \arrow[d, "\id \otimes \psi(\wedge^2 V_0)_{b,0}^*"] & \wedge^{2b-1}(V_1) \otimes \Sym^{b-1}(V_0)^* \arrow[l] \arrow[d, "\id \otimes \psi(\wedge^2V_0)_{b,0}^*"] & \wedge^{2b-1}(V_1) \otimes \Sym^{b-1}(V_0)^* \arrow[l] \arrow[d, "\id \otimes \psi(V_1)_{b,0}^*"]\\
\wedge^{2b-1}(\wedge^2V_0) \otimes \wedge^b(\wedge^2 V_0)^* \arrow[d] & \wedge^{2b-1}(V_1) \otimes \wedge^{b}(\wedge^2V_0)^* \arrow[l] \arrow[r] \arrow[d] & \wedge^{2b-1}V_1 \otimes \wedge^b V_1^* \arrow[d]\\
\wedge^{b-1}(\wedge^2 V_0) \arrow[d] & \wedge^{b-1}(V_1) \arrow[l] & \wedge^{b-1} V_1 \arrow[l] \arrow[dll, bend left = 10]\\
\mathbb{S}_{(b-1, 1^{b-1})}
\end{tikzcd}\]
Here the horizontal maps to the left are given by inclusions induced by $V_1 \subset \wedge^2(V_0)$ and the horizontal map to the right is a quotient map dual to such an inclusion. $g_{b-1}$ restricted to $\wedge^{2b-1}(V_1) \otimes \Sym^{b-1}(V_0)^*$ is the composition starting in the upper right and then traversing counterclockwise along the outside of the diagram. Likewise, $f(V_1)_{b-1} \circ h$ is the composition starting in the upper right and then traversing clockwise along the outside of the diagram. Each of the smaller square/triangles within the diagram is clearly commutative proving the commutativity of \eqref{containmentfactor}.
\end{proof}
\section{Special Hermite Isomorphism}\label{specHerm}
Crucial to the proof of \Cref{main1} at the end of \Cref{proofsec}, is the isomprhism $\psi_{b,0}(V_1)$. In the $SL_2(\mathbb{C})$-equavariant case when $V_0 = \Sym^b U$ and $V_1 = \Sym^{2b-2}U$, $\psi_{b,0}(\Sym^{2b-2}U)$ is a Hermite Reciprocity isomorphism 
\[\psi_{b,0}(V_1): \wedge^{b}(\Sym^{2b-2}U) \to \Sym^{b-1}(\Sym^b U) \otimes \det(\Sym^b U)\]
Throughout this paper, we have given three different characterizations of this map. The first is that it is a restriction of an $SL(V_0)$-equivariant projection
\[\wedge^b(\wedge^2V_0) \to \mathbb{S}_{(b, 1^b)}V_0\]
By \Cref{psidef} it factors as in the following diagram
\[\begin{tikzcd}\label{factorreform}
 \wedge^bV_1 \arrow[d] \arrow[ddr, bend left = 30, "\psi_{b,0}(V_1)"]\\
 \wedge^b(V_0 \otimes V_0) \arrow[d]\\
 \wedge^b V_0 \otimes \Sym^bV_0 & \arrow[l] \mathbb{S}_{b, 1^b} V_0
 \end{tikzcd}\]
Finally, we gave a cohomological interpretation in \eqref{cohomreform} by
\[\psi_{b,0}(V_0) = H^0(\wedge^{b}V_1\otimes \mathscr{O}_{\Pp(V_0)} \to \wedge^b K (b))\]
As mentioned in the introduction, our formulation of Hermite Reciprocity does not agree with the several formulations discussed in \cite{RaicuSam}. In the following two examples, we see this explicitly in the simplest case where the can differ, that is when $b = 3$ (in the case of $b = 2$, $\psi_{b,0}(V_1)$ is essentially just the identity map)

\begin{example}\label{spechermex}
Let $b = 3$ and $U$ have basis {x,1}. The map $\psi_{b, 0}(\Sym^{2b-2} U)$ is in this case an isomorphism
\[\psi_{b,0}(\Sym^4 U): \wedge^3(\Sym^4 U) \to \Sym^{2}(\Sym^3 U) \otimes \det(\Sym^3 U)\]
We will write the map $\psi_{3,0}(\Sym^{4}U)$ as a matrix using the basis for $U$ and the standard bases for exterior and symmetric powers. If we let $\{z_0, z_1, z_2 ,z_3\}$ correspond to the basis $\{x^3, x^2, x, 1\}$, then we can write the linear map of $S = \Sym^{\bullet}(\Sym^3 U)$-modules 
\[\Sym^{4}U \otimes S \to \Sym^3 \otimes S\]
as a matrix
\[\begin{array}
{c|ccccc}
 & x^4 & x^3 & x^2 & x & 1\\
 \hline
 x^3 & -z_1 & \frac{-1}{2}z_2& \frac{-1}{6}z_3 & 0 & 0\\
 x^2 & z_0 & 0 & \frac{-1}{2}z_2 & \frac{-1}{2}z_3 & 0\\
 x & 0 & \frac{1}{2}z_0 & \frac{1}{2} z_1 & 0 & -z_3\\
 1 & 0 & 0 & \frac{1}{6}z_0 & \frac{1}{2}z_1 & z_2
\end{array}\]

Taking the $3$x$3$ minors of this matrix yields a map
\begin{equation}\label{3x3}\wedge^3(\Sym^4 U) \to \wedge^3(Sym^3 U) \otimes \Sym^3(\Sym^3 U)\end{equation}
Using the basis for $U$ and the standard bases, we obtain a basis for
\[\wedge^{4}(\Sym^3 U) \otimes \Sym^2(\Sym^3U) \]
Mapping these basis vectors via the Koszul differential
\[\wedge^{4}(\Sym^3 U) \otimes \Sym^2(\Sym^3U) \to \wedge^3(Sym^3 U) \otimes \Sym^3(\Sym^3 U)\]
we obtain a basis for the subspace 
\[\wedge^{4}(\Sym^3U)\otimes \Sym^2(\Sym^3 U) = \mathbb{S}_{(3, 1^3)}(\Sym^3 U) \subset \wedge^3(\Sym^3U) \otimes \Sym^3 (\Sym^3 U)\]
The image of the map in \eqref{3x3}, can then be written using this basis. Ignoring the $1$-dimensional determinantal factor, $\wedge^4(\Sym^3 U)$, $\psi_{3,0}(\Sym^4 U)$ is given by the matrix
\[
\renewcommand{\arraystretch}{1.3} 
\begin{array}{c|@{\hskip -5pt}c@{\hskip -5pt}c@{\hskip -5pt}c@{\hskip -5pt}c@{\hskip -5pt}c@{\hskip -5pt}c@{\hskip -5pt}c@{\hskip -5pt}c@{\hskip -5pt}c@{\hskip -5pt}c}
  & \hspace{15pt}\rotatebox{60}{$x^4\wedge x^3 \wedge x^2$} 
  & \hspace{7pt}\rotatebox{60}{$x^4 \wedge x^3 \wedge x$} 
  & \hspace{7pt}\rotatebox{60}{$x^4\wedge x^3\wedge1$} 
  & \hspace{7pt}\rotatebox{60}{$x^4\wedge x^2 \wedge x$} 
  & \hspace{7pt}\rotatebox{60}{$x^4\wedge x^2 \wedge 1$} 
  & \hspace{7pt}\rotatebox{60}{$x^3\wedge x^2 \wedge x$} 
  & \hspace{7pt}\rotatebox{60}{$x^4\wedge x \wedge 1$} 
  & \hspace{7pt}\rotatebox{60}{$x^3\wedge x^2 \wedge 1$} 
  & \hspace{7pt}\rotatebox{60}{$x^3\wedge x \wedge 1$} 
  & \hspace{7pt}\rotatebox{60}{$x^2\wedge x \wedge 1$} \\
\hline
x^3\cdot x^3  & 1 & 0 & 0 & 0 & 0 & 0 & 0 & 0 & 0 & 0 \\
x^3 \cdot x^2  & 0 & 3 & 0 & 0 & 0 & 0 & 0 & 0 & 0 & 0 \\
x^3\cdot x  & 0 & 0 & 6 & 0 & 0 & 0 & 0 & 0 & 0 & 0 \\
x^2 \cdot x^2  & 0 & 0 & 0 & 3 & 0 & 0 & 0 & 0 & 0 & 0 \\
x^3 \cdot 1  & 0 & 0 & 0 & 0 & 2 & \frac{-1}{2} & 0 & 0 & 0 & 0 \\
x^2 \cdot x  & 0 & 0 & 0 & 0 & 6 & \frac{3}{2} & 0 & 0 & 0 & 0 \\
x^2 \cdot 1  & 0 & 0 & 0 & 0 & 0 & 0 & 6 & 0 & 0 & 0 \\
x\cdot x  & 0 & 0 & 0 & 0 & 0 & 0 & 0 & 3 & 0 & 0 \\
x\cdot 1 & 0 & 0 & 0 & 0 & 0 & 0 & 0 & 0 & 3 & 0 \\
1\cdot 1 & 0 & 0 & 0 & 0 & 0 & 0 & 0 & 0 & 0 & 1
\end{array}
\]
\end{example}

\begin{example}\label{oldhermex}
We can also write down the Hermite Reciprocity isomorphism found in \cite{RaicuSam} and elsewhere as a matrix in this case. \cite{Greens}[Section 3.4] gives an explicit description of an isomorphism, which specializes to a map
\[\Sym^2(D^3U) \to \wedge^3(\Sym^4 U).\]
Using again the basis $\{1,x\}$ for $U$, we can rewrite this map using the standard bases for symmetric and exterior powers, i.e. making the identification in characteristic 0 that
\[\Sym^2(D^3 U) \cong \Sym^2(\Sym^3 U)\]
Taking inverses, yields an isomorphism
\[\wedge^3(\Sym^4 U) \to \Sym^2 (\Sym^3 U)\], which is given by the matrix
\[
\renewcommand{\arraystretch}{1.3} 
\begin{array}{c|@{\hskip -5pt}c@{\hskip -5pt}c@{\hskip -5pt}c@{\hskip -5pt}c@{\hskip -5pt}c@{\hskip -5pt}c@{\hskip -5pt}c@{\hskip -5pt}c@{\hskip -5pt}c@{\hskip -5pt}c}
  & \hspace{15pt}\rotatebox{60}{$x^4\wedge x^3 \wedge x^2$} 
  & \hspace{7pt}\rotatebox{60}{$x^4 \wedge x^3 \wedge x$} 
  & \hspace{7pt}\rotatebox{60}{$x^4\wedge x^3\wedge1$} 
  & \hspace{7pt}\rotatebox{60}{$x^4\wedge x^2 \wedge x$} 
  & \hspace{7pt}\rotatebox{60}{$x^4\wedge x^2 \wedge 1$} 
  & \hspace{7pt}\rotatebox{60}{$x^3\wedge x^2 \wedge x$} 
  & \hspace{7pt}\rotatebox{60}{$x^4\wedge x \wedge 1$} 
  & \hspace{7pt}\rotatebox{60}{$x^3\wedge x^2 \wedge 1$} 
  & \hspace{7pt}\rotatebox{60}{$x^3\wedge x \wedge 1$} 
  & \hspace{7pt}\rotatebox{60}{$x^2\wedge x \wedge 1$} \\
\hline
x^3\cdot x^3  & 1 & 0 & 0 & 0 & 0 & 0 & 0 & 0 & 0 & 0 \\
x^3 \cdot x^2  & 0 & 3 & 0 & 0 & 0 & 0 & 0 & 0 & 0 & 0 \\
x^3\cdot x  & 0 & 0 & 3 & -3 & 0 & 0 & 0 & 0 & 0 & 0 \\
x^2 \cdot x^2  & 0 & 0 & 0 & 9 & 0 & 0 & 0 & 0 & 0 & 0 \\
x^3 \cdot 1  & 0 & 0 & 0 & 0 & -1 & 1 & 0 & 0 & 0 & 0 \\
x^2 \cdot x  & 0 & 0 & 0 & 0 & 9 & 0 & 0 & 0 & 0 & 0 \\
x^2 \cdot 1  & 0 & 0 & 0 & 0 & 0 & 0 & -3 & 3 & 0 & 0 \\
x\cdot x  & 0 & 0 & 0 & 0 & 0 & 0 & 0 & 9 & 0 & 0 \\
x\cdot 1 & 0 & 0 & 0 & 0 & 0 & 0 & 0 & 0 & 3 & 0 \\
1\cdot 1 & 0 & 0 & 0 & 0 & 0 & 0 & 0 & 0 & 0 & 1
\end{array}
\]
\end{example}

Comparing the matrices in \Cref{spechermex} and \Cref{oldhermex}, we see that they differ for the three $2$x$2$ matrices along the diagonal. This is possible because 
\[\wedge^3(\Sym^4 U) \cong \Sym^2{\Sym^3 U}\]
has $SL_2(\mathbb{C})$-equavariant decomposition into irreducibles
\[\Sym^6 U \oplus \Sym^2 U.\]
Thus, these representations have three weight spaces that are $2$-dimensional corresponding to these $2$x$2$ blocks.

\section*{Acknowledgements}
 The author would like to thank Claudiu Raicu for proposing this project and for numerous helpful discussions. The author further acknowledges the support of the National Science Foundation Grant DMS-2302341. Experiments in Macaulay2 \cite{GS} provided valuable insight for this project.

\bibliography{main}

\end{document}